\documentclass[12 pt, draftcls, onecolumn]{IEEEtran}
\usepackage{cite}

\ifCLASSINFOpdf
\else
\fi
\usepackage{amsmath}
\usepackage{algorithmic}

\usepackage{array}
\usepackage{fixltx2e}

\usepackage{stfloats}

\usepackage{algorithm}
\usepackage{algorithmic}
\usepackage{amsfonts}
\usepackage{graphicx}
\newtheorem{lemma}{Lemma}
\newtheorem{remark}{Remark}

\newtheorem{assumption}{Assumption}

\newtheorem{theorem}{Theorem}
\begin{document}
%
\title{Distributed Solver for Discrete-Time Lyapunov Equations Over Dynamic Networks with Linear Convergence Rate}
%
%
%

\author{Xia Jiang,~
        Xianlin~Zeng,~
        Jian Sun，~and~Jie Chen,~\IEEEmembership{Fellow,~IEEE}
}
\maketitle

\begin{abstract}
 This paper investigates the problem of solving discrete-time Lyapunov equations (DTLE) over a multi-agent system, where every agent has access to its local information and communicates with its neighbors. To obtain a solution to DTLE, a distributed algorithm with uncoordinated constant step sizes is proposed over time-varying topologies. The convergence properties and the range of constant step sizes of the proposed algorithm are analyzed. Moreover, a linear convergence rate is proved and the convergence performances over dynamic networks are verified by numerical simulations.
\end{abstract}

%
\IEEEpeerreviewmaketitle

\section{Introduction}
\par Distributed algorithms are vital for control and optimization over large-scale systems, such as sensor networks \cite{WuFu-56}, wireless communication networks \cite{shi_tradeoff} and smart grids \cite{OPF}. The research of distributed optimization algorithms has attracted considerable interest. Distributed optimization problems, which search for an optimal solution of a sum of individual objective functions with local or coupled constraints, have been widely investigated. A large number of effective discrete-time algorithms have been developed to solve distributed optimization problems in multi-agent networks, such as distributed gradient/subgradient methods \cite{H_average,consis_ste,IA96,LiangWang-98,NecoaraNedelcu-97}.
\par In recent years, many distributed optimization algorithms have been proposed for the computational problems of solving large-scale linear algebraic equations of the form $Ax=b$, where $A$ is a matrix and $x$, $b$ are vectors of appropriate dimensions \cite{Mou1,LiuMorse-53,Liu94,B_and,shi_net}. \cite{Mou1} provided a distributed computing algorithm and detailed discussions for solving linear algebraic equations, where every agent only knows a subset of the partitioned matrix $\left[\begin{matrix}A & b\end{matrix}\right]$ and exchanges the estimates of solution with its neighbors. Based on the assumption of the existence of at least one exact solution to linear algebraic equations, \cite{LiuMorse-53} proposed a distributed discrete-time algorithm with a linear convergence rate. If linear equations have no exact solutions, distributed algorithms to solve the least-squares solution are investigated \cite{Liu94,shi_net}. Linear matrix equations are more general than linear algebraic equations and the above distributed algorithms cannot be applied directly to solving linear matrix equations.

\par Linear matrix equations have wide applications in the design of modern complex systems. In particular, the discrete-time Lyapunov equation (DTLE), which is one of the most common linear matrix equations, is particularly important in the analysis of linear dynamic systems, such as stability analysis, controllability analysis and optimal control \cite{linear_theory,review_ly}. There is a large amount of literature on centralized numerical methods to solve Lyapunov equations \cite{Hammarling-91,BennerHeinkenschloss-79,WuSun-87,HaberVerhaegen-68}. With the growing data in systems and the expansion of the system scale, many excellent works have studied efficient methods for solving large-scale Lyapunov equations \cite{BennerHeinkenschloss-79,HaberVerhaegen-68}. However, because system information is distributed among separated equipments and computing capability of one agent is limited, the need to solve Lyapunov equations in a distributed way arises. Recently, there are several works studying distributed algorithms for computations of linear matrix equations \cite{zeng_opt,jia}. \cite{zeng_opt} has investigated the problem of solving general linear matrix equations and proposed distributed continuous-time algorithms over time-invariant undirected topologies. The distributed algorithm for solving DTLE proposed in \cite{jia} is a gradient-based method over time-invariant undirected graphs with diminishing step sizes, whose convergence rate is not fast enough.

\par This paper proposes a discrete-time distributed algorithm for solving DTLE of the form $AXA'-X+Q=0_{n \times n}$ over time-varying topologies, which has a linear convergence rate. The contributions are summarized below.  
\begin{itemize}
\item The paper studies the problem of solving discrete-time Lyapunov equation in a distributed way, which is vital in the analysis of linear dynamical systems. \cite{zeng_opt} proposed a distributed continuous-time algorithm over time-invariant undirected graphs, which cannot be applied over time-varying graphs. To our best knowledge, there is no work studying how to solve the discrete-time Lyapunov matrix equation over time-varying topologies in a distributed way.

\item This paper proposes a distributed discrete-time algorithm with uncoordinated constant step sizes over time-varying graphs. The proposed method overcomes the shortage of algorithms with diminishing step sizes that fail to achieve linear convergence rates \cite{RN68,RN69,jia}. \cite{consis_ste,ext} have developed algorithms with identical constant step sizes over multi-agent topologies. However, the methods can only achieve linear convergence rates under strongly convex assumptions. In this paper, the proposed method achieves a linear convergence rate without the assumption of strong convexity. What's more, the proposed algorithm has uncoordinated step sizes.

\item Using the theories of optimization and control, we provide sufficient conditions for the linear convergence rate of the proposed algorithm under time-varying topologies and the range of constant step sizes of the proposed algorithm.
\end{itemize}

 \par The remainder of the paper is organized as follows.
 \par The mathematical priliminaries and the problem description of solving discrete-time Lyapunov equation in a distributed manner are given in section \ref{problem_formulation}. A distributed discrete-time algorithm for the optimization problem and its convergent properties are provided in section \ref{solver_design}. We provide theoretical proofs for our main results in section \ref{proof_sec}. The convergence performances of the proposed algorithm are verified by numerical simulations in Section \ref{simulation} and the conclusion is made in section \ref{conclusion}.

\section{Priliminaries $\&$ Problem Description} \label{problem_formulation}
\subsection{Mathematical Priliminaries}
\par We denote $\mathbb{R}$ as the set of real numbers, $\mathbb{R}^n$ as the set of $n$-dimensional real column vectors, $\mathbb{R}^{n\times m}$ as the set of $n$-by-$m$ real matrices, $\mathbf{1}_m$ as an $m$-dimensional column vector with all elements being one and $E_n$ as the $n \times n$ identity matrix, respectively. We denote $A'$ as the transpose of matrix $A$, $0_{n\times q}$ as the $n \times q$ matrix with all elements of 0, $\lambda_{min}(A)$ as the minimum eigenvalue of the matrix $A$, $\lambda_{max}(A)$ as the maximum eigenvalue of the matrix $A$ and $A$ $\otimes$ $B$ as the Kronecker product of matrices $A$ and $B$. For a matrix $A$, the formula $A>0\quad (A \geq 0)$ denotes that $A$ is positive definite (positive semi-definite) and $A<0 \quad (A \leq 0)$ denotes that $A$ is negative definite (negative semi-definite). In addition, for a real vector $v$, $\left\|v\right\|$ is the Euclidean norm. For a real matrix $A$, $\left\|A \right\|_F$ denotes the Frobenius norm of the real matrix defined by $\left\|A\right\|_F=\sqrt{tr(A'A)}=\sqrt{\sum_{i,j}A_{ij}^2}$. Let $\left<\cdot,\cdot\right>_F$ be the Frobenius inner product of real matrices defined by $\left< B_1,B_2\right>_F=tr(B_1'B_2)=\sum_{i,j}(B_1)_{i,j}(B_2)_{i,j}$ with $B_1, B_2$ $\in$ $\mathbb{R}^{m\times n}$. If $f:\mathbb{R}^{m_1 \times m_2} \times \mathbb{R}^{n_1 \times n_2} \to \mathbb{R}$, $\nabla_{X}f(X,Y)=[\nabla_{X{ij}}f(X,Y)]$ $\in$ $\mathbb{R}^{m_1 \times m_2}$ denotes the partial gradient of function $f$ with respect to $X$. For a matrix $A \in \mathbb{R}^{n \times m}$, $vec(A)$ denotes the vectorization of $A$ which is an $mn \times 1$ column vector obtained by juxtaposing the consecutive rows of the matrix next to each other  and taking the transpose of the obtained long “multi-row”: $vec(A)=[a_{11},\cdots,a_{1m},a_{21},\cdots,a_{2m},\cdots,a_{n1},\cdots,a_{mn}]'$.

\subsection{Graph Theory}
\par The dynamic communication between $m$ agents at time $k$ is denoted by $\mathcal{G}_k(\nu,\varepsilon_k,\mathcal{A}_k)$, where $\nu=\{1,\dots,m\}$ is a finite nonempty node set, $\varepsilon_k$ $\subset  $ $\nu \times \nu$ is the edge set of time $k$ and $\mathcal{A}_k=[a_{ij,k}]$ $\in$ $\mathbb{R}^{m \times m}$ is the adjacency matrix at time $k$. For undirected graphs, all the adjacency matrices $\{\mathcal{A}_k\}$ are symmetric matrices such that $a_{ij,k}=a_{ji,k}$ and diagonal element $a_{ii,k}=0$. If an edge $(i,j) \in \varepsilon_k$, then node $j$ is called a neighbor of agent $i$ at time $k$ and $a_{ij,k}>0$, otherwise, $a_{ij,k}=0$. Let $\mathcal{N}_{i,k}$ denotes the set of neighbor nodes that connect to node $i$ at time $k$, i.e., $\mathcal{N}_{i,k}=\{j:(i,j)\in \varepsilon_k\}$. The Laplacian matrix $\mathcal{L}_k$ at time $k$ is defined as $l_{ij,k}=-a_{ij,k}$, for $i \neq j$ and $l_{ii,k}=\sum_{j \in \mathcal{N}_{i,k}} a_{ij,k}$, which ensures that $\sum_{j=1}^m l_{ij,k}=0$. What's more, $\mathcal{L}_k =\mathcal{L}_k'$ for every undirected graph $\mathcal{G}_k$. The dynamic undirected graphs are assumed to be simple which implies that there is no repeated edge or self-loop.
\subsection{DTLE Problem Description}
\par The discrete-time Lyapunov equation (DTLE) for stability/controllability \cite{linear_theory} is
\begin{align}\label{target1}
AXA'-X+Q=0_{n \times n},
\end{align}
where $X$, $A$, $Q \in \mathbb{R}^{n\times n}$ and $X$ is the unknown variable. In addition, $X$ and $Q$ are symmetric matrices. If the Lyapunov equation (\ref{target1}) is used to analyze one system's stability, $Q$ is a symmetric positive definite matrix.
\par The matrices $A$ and $Q$ are decomposed as follows:
\begin{align} \label{matrix_decomposition}
A=\left[ \begin{matrix} A_{r1}\\ \cdot \\ \cdot \\ \cdot \\A_{rm} \end{matrix} \right] \in \mathbb{R}^{n \times n}, \quad Q=\left[ \begin{matrix} Q_{l1} & \cdot & \cdot & \cdot & Q_{lm} \end{matrix} \right] \in \mathbb{R}^{n \times n},
\end{align}
where $A_{ri}\in \mathbb{R}^{n_{i} \times n}$, $Q_{li} \in \mathbb{R}^{n \times n_{i}}$, $\sum_{i=1}^m n_{i}=n$.
\par In this structure, every agent $i$ only knows $A_{ri}$, $Q_{li}$, and communicates with its neighbors $\{j: j\in \mathcal{N}_i\}$ to obtain a same solution to (\ref{target1}). To distinguish the row-blocks or column-blocks of a matrix, we use subscript $ri$ to denote its $i$th row-block and subscript $li$ to denote its $i$th column-block.

\par The objective of this paper is to propose a fully distributed discrete-time algorithm over time-varying undirected topologies to solve DTLE (\ref{target1}) with a linear convergence rate. We need the following assumptions.
\begin{assumption}\label{assumption1}
Equation \eqref{target1} has at least one solution.
\end{assumption}
\begin{assumption}\label{assumption3}
Graph $\mathcal{G}_k$ is undirected for all $k \geq 0$ and the adjacency matrix $\mathcal{A}_k$ of graph $\mathcal{G}_k$ is doubly stochastic.
\end{assumption}
\begin{remark} Assumption \ref{assumption1} makes the problem well-defined, which also implies that the solution set is non-empty. Assumption \ref{assumption3} is common for the design of distributed algorithms over multi-agent topologies.
\end{remark}

\section{Distributed Solver Design}\label{solver_design}

In this section, we reformulate the problem into a separable one and propose a distributed discrete-time algorithm with a linear convergence rate.
\subsection{Problem Reformulation}\label{pro_trans}
To deal with the coupling in $AXA'$, we introduce an auxiliary variable $Y \in \mathbb{R}^{n \times n}$ to transform the term into a separable structure. Specifically, equation (\ref{target1}) is equivalent to
\begin{align*}
Y&=AX,\\
YA'&=XE_n-Q,
\end{align*}
where $E_n$ is an $n \times n$ identity matrix.
\par Define $X_{i} \in \mathbb{R}^{n \times n}$ and $Y_{i} \in \mathbb{R}^{n \times n}$ as estimates of $X$ and $Y$ of agent $i$ $\in$ $\{1,...,m\}$ respectively. With the matrix decomposition (\ref{matrix_decomposition}), $Y=AX$ can be transformed into $Y_{i}^{ri}=A_{ri}X_i$ and $Y_i=[Y_{r1}', \cdots, Y_{rm}']'$; By multiplying $X$ with an identity matrix $E_n$, $YA'=XE_n-Q$ can be transformed into $Y_iA_{ri}'=X_{i}E_{li}-Q_{li}$ and $E_n=[E_{l1},\cdots,E_{lm}]$. With the above decompositions, we decentralize the DTLE (\ref{target1}) into
\begin{align}\label{distributed_equations}
\begin{split}
Y_{i}^{ri}=&A_{ri}X_{i}, \\
Y_{i}A_{ri}'=&X_{i}E_{li}-Q_{li},\\
X_{i}=&X_{j}, \\
Y_{i}=&Y_{j},
\end{split}
\end{align}
for all $i$, $j$ $\in$ $\{1,...,m\}$.
\par Then, we further transform the distributed computation of (\ref{distributed_equations}) into a distributed optimization problem.
\begin{align}\label{opti_problem_origin}
  &Min._{(\hat X,\hat Y)} \quad F(\hat X,\hat Y)=\sum_{i=1}^{m}f_{i}(X_i,Y_i), \notag\\
  &s.t. \quad X_i=X_j, \quad Y_i=Y_j, \quad i,j \in \{1,\cdots, m\},
\end{align}
where $f_{i}(X_i,Y_i)=\frac{1}{2}\|Y_i^{ri}-A_{ri}X_i \|_F ^2+\frac{1}{2} \|Y_{i}A_{ri}'-X_iE_{li}+Q_{li} \|_F ^2$, $\hat Y = [Y_1^{\rm '},\ldots,Y_n^{\rm '}]^{\rm '}\in\mathbb R^{mn\times n}$, $\hat X = [X_1^{\rm '},\ldots,X_n^{\rm '}]^{\rm '}\in\mathbb R^{mn\times n}$, agent $i$ only knows the information of $A_{ri}$, $E_{li}$ and $Q_{li}$. Throughout this paper, we use $(\hat{X}^*,\hat{Y}^*)$ to denote an optimal solution to problem (\ref{opti_problem_origin}).

\begin{remark}
Because the optimal function value of (\ref{opti_problem_origin}) $F^*=0$, it implies that four equations in (\ref{distributed_equations}) are satisfied. Then an optimal solution of problem (\ref{opti_problem_origin}) is $\hat{X}^*$ if and only if $\hat{X}^*=\mathbf{1}_m \otimes X^*$ and $X^*$ is a solution of (\ref{target1}).
\end{remark}

\subsection{Distributed algorithm}\label{algo_propose_sec}
\par From problem reformulation, every agent $i$ has the information state $X_{i,k}\in \mathbb{R}^{n \times n}$ and the auxiliary variable $Y_{i,k}\in \mathbb{R}^{n \times n}$, which are the estimates of the solutions of the original optimization problem (\ref{opti_problem_origin}) at time $k$.

\par Agent $i$ updates its estimates $X_{i,k+1}$ according to the following rules. For all $k\geq 0$ and all $i \in \{1,...,m\}$:
\begin{small}
\begin{align}\label{algo_full}
&X_{i,k+1}=X_{i,k}-\alpha_i d_{X_{i,k}}\!-\!\frac{\alpha_i}{2}\!\sum_{j\in \mathcal{N}_{i,k}}\! a_{ij,k}(X_{i,k}\!-\!X_{j,k}), \notag\\
&Y_{i,k+1}=Y_{i,k}\!-\!\alpha_i d_{Y_{i,k}}\!-\!\frac{\alpha_i}{2}\!\sum_{j\in \mathcal{N}_{i,k}}\! a_{ij,k}(Y_{i,k}\!-\!Y_{j,k}), \notag\\
&d_{X_{i,k}}\!=\!-\!A_{ri}'(Y_{i,k}^{ri}\!-\!A_{ri}X_{i,k})\!-\!(Y_{i,k}A_{ri}'\!-\!X_{i,k}E_{li}\!+\!Q_{li})E_{li}', \notag \\
&d_{Y_{i,k}}\!=\!E_{li}(Y_{i,k}^{ri}\!-\!A_{ri}X_{i,k})\!+\!(Y_{i,k}A_{ri}'\!-\!X_{i,k}E_{li}\!+\!Q_{li})A_{ri},
\end{align}
\end{small}
where $a_{ij,k}$ is the weight of time $k$ between agent $i$ and $j$, $d_{X_{i,k}}=\nabla_{X_{i,k}}f_i(X_{i,k},Y_{i,k})$, which is a gradient of the objective function $f_i(X_{i},Y_{i})$ with respect to $X_i$ at time $k$ and similarly, $d_{Y_{i,k}}=\nabla_{Y_{i,k}}f_i(X_{i,k},Y_{i,k})$.
\par For ease of presentation, we introduce an assumption to describe the range of step sizes.
\begin{assumption}\label{xi_defi}
In the proposed algorithm (\ref{algo_full}), the step size of every agent $i$, $\alpha_i$ satisfies that min$\{1,1/\xi_i\}>\alpha_i>0$, where $\xi_i$ satisfies the following inequality
\begin{align*}
\|\!-\!A_{ri}'T_{i1}\!-\!T_{i2}E_{li}'\|^2_{F}&+\|E_{li}T_{i1}\!+\!T_{i2}A_{ri}\|^2_{F}\\
&\leq \xi_i\|T_{i1}\|^2_{F}\!+\!\xi_i\|T_{i2}\|^2_{F},
\end{align*}
for all $T_{i1}\in\mathbb R^{n_i\times n}$, $T_{i2}\in\mathbb R^{n\times n_i}$.
\end{assumption}
\begin{remark}
Because of the special form of the above inequality, it is easy to find a proper variable $\xi_i$. For instance, by transformation, every $\xi_i$ satisfying $\xi_i \geq 2 (\|A_{ri}'\|^2+\|E_{li}\|^2)$ also satisfies the above inequality. So, we can find a proper step size $\alpha_i$ satisfying Assumption \ref{xi_defi} easily.
\end{remark}
For the simplicity of analysis, we define $Z_{i,k} = \left[\begin{matrix}Y_{i,k}\\X_{i,k}\end{matrix}\right]$, $\hat{Z}_k=\left[\begin{matrix}Z_{1,k}\\\vdots\\Z_{m,k}\end{matrix}\right]$,  and $f_i(Z_i)=f_{i}(X_i,Y_i)$ without causing ambiguity.
Rewrite the algorithm (\ref{algo_full}) compactly in terms of the matrix $Z_{i,k}$ as
\begin{align*}
 Z_{i,k+1}&=Z_{i,k}\!-\!\alpha_i d_{Z_{i,k}}\!-\!\frac{\alpha_i}{2}\!\sum_{j\in \mathcal{N}_{i,k}}\! a_{ij,k}(Z_{i,k}-Z_{j,k})\\
&=(1-\frac{\alpha_i}{2}l_{ii,k})Z_{i,k}\!-\!\frac{\alpha_i}{2}\!\sum_{j\in \mathcal{N}_{i,k}} \!l_{ij,k}Z_{j,k}\!-\!\alpha_i d_{Z_{i,k}},
 \end{align*}
where $i$ $\in$ $\{1,...,m\}$, $k\geq 0$ and $d_{Z_{i,k}}=\nabla_{Z_{i,k}} f_i(Z_{i,k})$.

\par We introduce a matrix $W_k=[w_{ij,k}]\in \mathbb{R}^{m \times m}$ at time $k$ and let $w_{ii,k}=1-\frac{\alpha_i}{2}l_{ii,k}$ and $w_{ij,k}=-\frac{\alpha_i}{2}l_{ij,k}$. Then the above equation can be further written as
\begin{align}\label{algo_com}
Z_{i,k+1}=\sum_{j=1}^m w_{ij,k}Z_{j,k}-\alpha_i \nabla f_i(Z_{i,k}),
\end{align}
where $\sum_{j=1}^m w_{ij,k}=1$ and the matrix $W_k\in \mathbb{R}^{m \times m}$ is doubly stochastic if Assumption \ref{assumption3} holds.

\begin{remark}
There are a lot of excellent works studying distributed optimization algorithms with diminishing/constant step sizes in detail. In comparison, the algorithm proposed in this paper has the following differences and advantages.
\begin{itemize}
\item The step sizes of many algorithms are diminishing, hindering the possibility of linear convergence rates \cite{RN68,jia}. The proposed algorithm has constant step sizes and a linear convergence rate.
\item Some studies with constant step sizes have achieved linear convergence rates for strongly convex functions \cite{consis_ste,ext}. In addition, these methods require an identical step size for all agents. Nevertheless, the constant step sizes in the proposed algorithm are not necessarily identical and the method achieves a linear convergence rate without the strongly convex assumptions, making the algorithm more practical.
\item In \cite{au_gradient,zeng_opt}, they studied distributed algorithms over time-invariant undirected graphs. In comparison, the proposed algorithm of this paper can be applied to time-varying undirected topologies. Furthermore, the proposed algorithm uses no dual variable, making it more straightforward than primal-dual methods in \cite{zeng_opt}.
\end{itemize}

\end{remark}
\subsection{Convergence properties}\label{convergent}

In this subsection, convergence properties and the convergence rate of the proposed algorithm are presented.
\par To derive convergence properties of the proposed algorithm (\ref{algo_full}), we need some common assumptions.
\begin{assumption}\label{assumption2}
The graph sequence \{$\mathcal{G}_k$\} is uniformly connected. That is, there exists an integer $B\geq 1$ such that agent $i$ sends its information to all other agents at least once every $B$ consecutive time slots.
\end{assumption}
\begin{assumption}\label{assumption_w3}
There exists a scalar $\eta$ with $0<\eta<1$ such that for all $i \in \{1, \cdots m\}$, $w_{ii,k}\geq \eta$ for all $k\geq 0$ and $w_{ij,k} \geq \eta$ if agent $j$ communicates directly with agent $i$ at time $k$. Otherwise, $w_{ij,k}=0$.
\end{assumption}

\par Then the following theorem shows convergence results of the the proposed method.
\begin{theorem} \label{converge_theo}
Suppose Assumptions \ref{assumption1}, \ref{assumption3}, \ref{xi_defi} $-$ \ref{assumption_w3} hold.
 \begin{itemize}
 \item [(i)] The variable state of every agent $i$ generated by algorithm (\ref{algo_full}) converges to one same solution, which means $\hat{Z}_k\to \hat{Z}^*$ as $k \to \infty$ and $\hat{Z}^*$ is an optimal solution to problem (\ref{opti_problem_origin}).
 \item [(ii)] $\sum_{i=1}^m (1-\alpha_i \xi_i) f_i(\bar{Z}_{i,k})=O(\frac{1}{k})$, where $\bar{Z}_{i,k}$ is time average state variable of agent $i$, $\bar{Z}_{i,k}=\frac{\sum_{l=1}^k Z_{i,l}}{k}$.
 \end{itemize}
\end{theorem}

Theorem \ref{converge_theo} shows that over uniformly strongly connected graphs, all agents' variables converge to one same optimal solution. Next, we give sufficient conditions that the proposed algorithm has a linear convergence rate.

\begin{assumption}\label{connected_ass}
The topological sequence has finite graphs and every graph is connected.
\end{assumption}
\begin{remark}
In Assumption \ref{connected_ass}, the graph at every iteration $k$ needs to be connected and the number of graphs in topological sequence is finite, which is stronger than the conditions of Assumptions \ref{assumption2} and \ref{assumption_w3}. One topological sequence satisfying Assumption \ref{connected_ass} must satisfy Assumptions \ref{assumption2} and \ref{assumption_w3}, but the converse is not true.
\end{remark}
\par The linear convergence rate of the proposed algorithm is presented in the following theorem.

\begin{theorem}\label{conv_exp}
Suppose Assumptions \ref{assumption1}, \ref{assumption3}, \ref{xi_defi} and \ref{connected_ass} hold. The variable state $\hat{Z}_k$ generated by algorithm (\ref{algo_com}) converges to an optimal solution $\hat{Z}^*$ of problem (\ref{opti_problem_origin}) at a linear convergence rate.
\end{theorem}

\section{Theoretical Proof}\label{proof_sec}
In this section, we present theoretical proofs for the main results on convergence properties and convergence rate of algorithm (\ref{algo_full}).
\subsection{Proof of Theorem \ref{converge_theo}}
Before we give the proof of Theorem \ref{converge_theo}, we introduce a few lemmas.
\begin{lemma}  \label{eigenvalue_co}

Weyl's inequality: For matrix formula $M=H+P$, if $M,H$ and $P$ are all $n$ by $n$ Hermitian matrices (symmetric matrices for real matrices), the $i$th eigenvalue of $M$ satisfies that $\lambda_i(H) +\lambda_{min} (P)\leq \lambda_i(M)\leq \lambda_i(H)+\lambda_{max}(P)$.

\end{lemma}
\begin{lemma}\label{bound_lemma}
Suppose Assumptions \ref{assumption1}, \ref{assumption3} and \ref{xi_defi} hold.
\begin{itemize}
  \item [(i)] the sequence $\{\hat{Z}_k\}$ generated by algorithm (\ref{algo_com}) is bounded and every equilibrium of the sequence is Lyapunov stable.
  \item [(ii)] $\lim_{k\to \infty} f_i(Z_{i,k})=0$ for all agent $i$ and $\sum_{i=1}^m (1-\alpha_i \xi_i) f_i(\bar{Z}_{i,k})=O(\frac{1}{k})$, where $\bar{Z}_{i,k}$ is time average state variable of agent $i$, $\bar{Z}_{i,k}=\frac{\sum_{l=1}^k Z_{i,l}}{k}$.
  \item [(iii)] $\lim_{k\to \infty} e_{i,k}=0$ where $e_{i,k}=\alpha_i \nabla_{Z_{i,k}}f_i(Z_{i,k})$ for all agent $i$.
\end{itemize}
\end{lemma}
\textbf{Proof:}	
(i) Let $\hat{Z}^*$ be a solution to problem \eqref{opti_problem_origin}, which is also an equilibrium of the algorithm (\ref{algo_com}). Define $$ D(k) =\| \Lambda^{-\frac{1}{2}} ( \hat{Z}_{k+1}-\hat{Z}^*)\|^2_{F}- \| \Lambda^{-\frac{1}{2}} ( \hat{Z}_k-\hat{Z}^*)\|^2_{F},$$
 where $\Lambda=diag\{ \alpha_1\otimes E_{2n}, \cdots, \alpha_m \otimes E_{2n}\}$. Recall that ${L}_k=\mathcal{L}_k \otimes E_{2n}$. Since $\hat{Z}_{k+1} = \hat{Z}_k-\frac{1}{2} \Lambda L_k\hat{Z}_k-\Lambda \nabla_{\hat{Z}_k} F(\hat{Z}_k)$ by algorithm (\ref{algo_com}), we plug $\hat{Z}_{k+1}$ in $D(k)$ and get the following equation
\begin{align*}
 D(k) =& 2\big \langle \!-\!\frac{1}{2}{L}_k\hat{Z}_k\!-\! \nabla_{\hat{Z}_k} F(\hat{Z}_k), \hat{Z}_k\!-\!\hat{Z}^*\big\rangle_{F}\\
 &+\big\|\Lambda^{-\frac{1}{2}}(\frac{1}{2}\Lambda {L}_k\hat{Z}_k+\Lambda \nabla_{\hat{Z}_k} F(\hat{Z}_k))\big\|_F^2\\
 = & -\big \langle \hat{Z}_k-\hat{Z}^*, {L}_k\hat{Z}_k\big\rangle_{F} - 4F(\hat{Z}_k)\\
 &+\big\|\frac{1}{2}\Lambda^{\frac{1}{2}} {L}_k\hat{Z}_k+\Lambda^{\frac{1}{2}} \nabla_{\hat{Z}_k} F(\hat{Z}_k))\big\|_F^2.
\end{align*}
What's more, by triangle inequality, the last term satisfies that $\big\|\frac{1}{2}\Lambda^{\frac{1}{2}}{L}_k\hat{Z}_k+\Lambda^{\frac{1}{2}}\nabla_{\hat{Z}_k} F(\hat{Z}_k)\big\|^2_{F} \leq$ $ \frac{1}{2}\|\Lambda^{\frac{1}{2}}{L}_k\hat{Z}_k\|_{F}^2+2\|\Lambda^{\frac{1}{2}}\nabla_{\hat{Z}_k} F(\hat{Z}_k)\|^2_{F}$. We can obtain

\begin{align}\label{Dk_eq}
D(k) &\leq -\big \langle \hat{Z}_k-\hat{Z}^*, {L}_k\hat{Z}_k\big\rangle_{F}- 4F(\hat{Z}_k) \notag\\
 &+\frac{1}{2}\| \Lambda^{\frac{1}{2}}{L}_k\hat{Z}_k\|_{F}^2+2\| \Lambda^{\frac{1}{2}}\nabla_{\hat{Z}_k} F(\hat{Z}_k)\|^2_{F}\notag\\
 =& -\big \langle \hat{Z}_k-\hat{Z}^*, {L}_k(\hat{Z}_k-\hat{Z}^*)\big\rangle_{F}\notag\\
 &+\frac{1}{2}\big \langle\hat{Z}_k-\hat{Z}^*,L'_k\Lambda {L}_k(\hat{Z}_k-\hat{Z}^* )\big\rangle_{F}\notag\\
 &- 4F(\hat{Z}_k)+2\| \Lambda^{\frac{1}{2}}\nabla_{\hat{Z}_k} F(\hat{Z}_k)\|^2_{F}\notag\\
  =&\big\langle \hat{Z}_{k}-\hat{Z}^*, (\frac{1}{2}L'_k\Lambda {L}_k- {L}_k) (\hat{Z}_k-\hat{Z}^*) \big\rangle_{F}\notag\\
  &- 4F(\hat{Z}_k)+2\| \Lambda^{\frac{1}{2}}\nabla_{\hat{Z}_k} F(\hat{Z}_k)\|^2_{F}.
\end{align}
\par Firstly, as we know, $L_k$ is symmetric positive semi-definite for undirected graphs. So, there is a symmetic matrix $L_k^{\frac{1}{2}}$, which satisfies $L_k^{\frac{1}{2}}L_k^{\frac{1}{2}}=L_k$ \cite{D_nonopti}. $L_k^{\frac{1}{2}}$ is unique and called the square root of $L_k$. Then, we have
\begin{align}\label{L_max}
&L_k-\frac{1}{2}L'_k\Lambda {L}_k \notag\\
=&L_k^{\frac{1}{2}}(E-\frac{L_k^{\frac{1}{2}}\Lambda L_k^{\frac{1}{2}}}{2})L_k^{\frac{1}{2}}.
\end{align}
\par For the diagonal step-size matrix $\Lambda$, since that $(E-\Lambda) \geq 0$, then $L_k^{\frac{1}{2}}L_k^{\frac{1}{2}}-L_k^{\frac{1}{2}} \Lambda L_k^{\frac{1}{2}}\geq 0$. Furthermore, $(E-\frac{L_k^{\frac{1}{2}}\Lambda L_k^{\frac{1}{2}}}{2})-(E-\frac{L_k}{2})\geq 0$. The absolute value of eigenvalues of the doubly stochastic matrix is no larger than $1$ and ${L}_k=(E-\mathcal{A}_k)\otimes E_{2n}$. By Lemma \ref{eigenvalue_co}, the maximum eigenvalue of the Laplacian matrix ${L}_k$ is $2$, then $\lambda_{max}(\frac{1}{2}{L}_k) = 1$. It is obvious that $E-\frac{L_k}{2}\geq 0$ by Lemma \ref{eigenvalue_co}.
 \par What's more, because $(E-\frac{L_k^{\frac{1}{2}}\Lambda L_k^{\frac{1}{2}}}{2})-(E-\frac{L_k}{2})\geq 0$ and $E-\frac{L}{2}\geq 0$, we can also obtain that $(E-\frac{L_k^{\frac{1}{2}}\Lambda L_k^{\frac{1}{2}}}{2})\geq 0$ by Lemma \ref{eigenvalue_co}. Therefore, the matrix formula (\ref{L_max}) is positive semi-definite and $(\frac{1}{2}L'_k\Lambda {L}_k- {L}_k)\leq 0$, which implies that the first term of equation (\ref{Dk_eq}) is non-positive.

\par Secondly, recall the definition of step size limitation $\xi_i$ in Assumption \ref{xi_defi}. Assume that $T_{i1}=Y_{i,k}^{ri}\!-\!A_{ri}X_{i,k}$ and $T_{i2}=Y_{i,k}A_{ri}'-X_{i,k}E_{li}+Q_{li}$. Then, by the definition of $f_i(X_i,Y_i)$ and Assumption \ref{xi_defi}, we have
\begin{align}\label{Lamb_F}
&\|\!-\!A_{ri}'T_{i1}\!-\!T_{i2}E_{li}'\|^2_{F}\!+\! \|E_{li}T_{i1}\!+\!T_{i2}A_{ri}\|^2_{F}\! \leq\! 2 \xi_i f_i(X_i,\!Y_i),\notag\\
&\| \nabla_{Z_{i,k}} f_i(Z_{i,k})\|_F^2 \leq 2 \xi_i f_i(Z_{i,k}),\notag\\
&\| \nabla_{\hat{Z}_k} F(\hat{Z}_k)\|_F^2 \leq 2 \sum_{i=1}^m \xi_i f_i(Z_{i,k}),\notag\\
&\|\Lambda^{\frac{1}{2}} \nabla_{\hat{Z}_k} F(\hat{Z}_k)\|_F^2 \leq 2 \sum_{i=1}^m \alpha_i \xi_i f_i(Z_{i,k})\leq 2F(\hat{Z}_k).
\end{align}
Therefore, $2\|\Lambda^{\frac{1}{2}}  \nabla_{\hat{Z}_k} F(\hat{Z}_k)\|^2_{F}\leq 4 F(\hat{Z}_k)$.

In summary, due to the choice of $\alpha_i$, it follows that $\left<\hat{Z}_k-\hat{Z}^*,(\frac{1}{2}L'_k\Lambda {L}_k- {L}_k)(\hat{Z}_k-\hat{Z}^*)\right>_F\leq 0$ and $- 4F(\hat{Z}_k)+2\|\Lambda^{\frac{1}{2}} \nabla_{\hat{Z}_k} F(\hat{Z}_k)\|^2_{F}\leq 0$. Therefore, $D(k)\leq 0$, which implies that the sequence $\{\hat{Z}_k\}$ is bounded and $\hat{Z}^*$ is a Lyapunov stable equilibrium of (\ref{opti_problem_origin}).

(ii) From $\|\Lambda^{\frac{1}{2}} \nabla_{\hat{Z}_k} F(\hat{Z}_k)\|_F^2 \leq 2 \sum_{i=1}^m \alpha_i \xi_i f_i(Z_{i,k})$ in (\ref{Lamb_F}), we deduce that $2\|\Lambda^{\frac{1}{2}} \nabla_{\hat{Z}_k} F(\hat{Z}_k)\|^2_{F}- 4F(\hat{Z}_k) \leq -4 \sum_{i=1}^m (1-\alpha_i \xi_i) f_i(Z_{i,k})\leq 0$. Denote $4 \sum_{i=1}^m (1-\alpha_i \xi_i) f_i(Z_{i,k})$ as $a_k$ such that  $a_k \geq 0$. Summing (\ref{Dk_eq}) over $k$, we have
\begin{align}\label{as_equation}
\|\Lambda^{-\frac{1}{2}}(\hat{Z}_{k+1}\!-\!\hat{Z}^*)\|_F^2\!-\!\|\Lambda^{-\frac{1}{2}}(\hat{Z}_0\!-\!\hat{Z}^*)\|_F^2 \leq \!-\!\sum_{s=1}^k\! a_s.
\end{align}
Since the variable sequence is bounded in (i), then we have $\sum_{s=1}^{\infty} a_s < \infty$ as $k \to \infty$, which implies $\lim_{k \to \infty}a_k=0$. By the definition of $a_k$, we have that $\lim_{k\to \infty} f_i(Z_{i,k})=0$ for all agent $i$.
\par Because the function $f_i$ is convex, it follows from Jensen's inequality \cite{Boyd_convex} that $4 k \sum_{i=1}^m (1-\alpha_i \xi_i) f_i(\bar{Z}_{i,k})\leq \sum_{s=1}^k a_s$. By equation (\ref{as_equation}),
\begin{align*}
\sum_{i=1}^m (1-\alpha_i \xi_i) f_i(\bar{Z}_{i,k})\leq \frac{\left\|\Lambda^{-\frac{1}{2}}(\hat{Z}_0-\hat{Z}^*)\right\|_F^2}{4k},
\end{align*}
and hence, $\sum_{i=1}^m (1-\alpha_i \xi_i) f_i(\bar{Z}_{i,k})=O(\frac{1}{k})$.
\par (iii) By part (ii), $\lim_{k\to \infty} f_i(Z_{i,k})=0$ for all agent $i$. Since $e_{i,k}=\alpha_i \nabla_{Z_{i,k}}f_i(Z_{i,k})$, it is straightforward that $\lim_{k\to \infty} e_{i,k}=0$ for all agent $i$.
$\hfill\square$

Recall that $Z_{i,k+1}=\sum_{j=1}^m w_{ij,k}Z_{j,k}-e_{i,k}$, where $e_{i,k}=\alpha_i \nabla_{Z_{i,k}}f_i(Z_{i,k})$. Define a transition matrix $\Phi(k,s)$ from time $s$ to time $k$ such that
\begin{align}
\Phi(k,s)=W_s W_{s+1}\cdots W_{k-1} W_k,
\end{align}
where $\Phi(k,k)=W_k$ for all $k$. It follows from the transition matrix and the proposed algorithm (\ref{algo_com}) that the relation between $Z_{i,k+1}$ and the estimates $Z_{1,s},\cdots, Z_{m,s}$ at time $s \leq k$ is given by
\begin{align}\label{trans_equa}
Z_{i,k+1}&=\sum_{j=1}^m \Phi(k,s)_{ij}Z_{j,s}\notag\\
&+\sum_{r=s+1}^k(\sum_{j=1}^m \Phi(k,s)_{ij} e_{j,r-1}) +e_{i,k}.
\end{align}

We define an auxiliary sequence $\{H_k\}$, where $H_k$ is given by $H_k=\frac{1}{m} \sum_{i=1}^m \sum_{j=1}^m w_{ij,k}Z_{j,k}$. Since the matrix $W(k)$ is doubly stochastic, it follows that:
\begin{align*}
H_k=\frac{1}{m} \sum_{j=1}^m Z_{j,k}.
\end{align*}
\par Furthermore, by the doubly stochasticity of the transition matrix, it follows from the relations equation (\ref{trans_equa}) that
\begin{align}\label{average_equa}
H_k=\frac{1}{m} \sum_{j=1}^m Z_{j,s}+\frac{1}{m} \sum_{r=s+1}^k(\sum_{j=1}^m e_{j,r-1}), \quad \forall k\geq s.
\end{align}
\par The next lemma shows that the limiting behavior of the agent's estimate $Z_{i,k}$ is the same as the limiting behavior of $H_k$ as $k \to \infty$.
\begin{lemma} \label{common_lemma}
Let Assumptions \ref{assumption1}, \ref{assumption3}, \ref{assumption2} and \ref{assumption_w3} hold. For all agent $i$, $\lim_{k \to \infty} \left\| Z_{i,k}-H_k \right\|_F=0$.
\end{lemma}
\textbf{Proof:} Because the matrix $W(k)$ is a doubly stochastic matrix, the transition matrix $\Phi(k,s)$ is also a doubly stochastic matrix. By Proposition $1$ in \cite{H_average}, the entries $[\Phi(k,s)]_{ij}$ of the transition matrix converge to $\frac{1}{m}$ as $k \to \infty$ with a linear rate for all $i,j \in \{ 1 \cdots m\}$,
\begin{align}\label{phi_lim}
\left| [\Phi(k,s)]_{ij}-\frac{1}{m} \right| \leq 2 \frac{1+ \eta^{-B_0}}{1-\eta^{B_0}}(1- \eta^{B_0})^{\frac{k-s}{B_0}},
\end{align}
for all $k \geq s$, where $\eta$ is the lower bound of Assumption \ref{assumption_w3}, $B_0=(m-1)B$ and $B$ is the intercommunication interval bound of Assumption \ref{assumption2}. By the Lemma \ref{bound_lemma} (iii), we have $\lim_{k\to \infty} e_{i,k}=0$ for all agent $i$. By the relationship (\ref{phi_lim}) and Lemma $4$ in \cite{H_average}, we get that $\lim_{k \to \infty} \left\| Z_{i,k}-H_k \right\|_F=0$.
$\hfill\square$

\par We have showed that the sequence $\{\hat{Z}_k\}$ is bounded and any optimal solution in the optimal solution set is Lyapunov stable. What's more, every agent's state $Z_{i,k}$ converges to one same individual average $H_k$ as $k \to \infty$. Next, we will prove that the sequence $\{\hat{Z}_k\}$ converges to one same optimal solution in the optimal solution set. Before that, we need a related lemma.
\begin{lemma}\label{lmm-conv}
Let the sequence $\{\hat{Z}_k\}$ be generated by  algorithm \eqref{algo_full}. Define $\Omega(\hat{Z}(\cdot))$ as the positive limit set. If $\hat Z\in\Omega(\hat{Z}(\cdot))$ is a Lyapunov stable equilibrium point of algorithm \eqref{algo_full}, then $\hat Z=\lim_{k \to \infty}\hat Z_k$ and $\Omega(\hat Z(\cdot)) = \{\hat Z\}$.
\end{lemma}
\textbf{Proof:}
The proof is similar to that of Proposition 3.1 in \cite{semi_sta} and is omitted.
$\hfill\square$
\par Then we are ready to give the proof of Theorem \ref{converge_theo}.
\par {\bf Proof of Theorem \ref{converge_theo}:}
From Lemmas \ref{bound_lemma} (ii), (iii) and \ref{common_lemma}, $\{\hat{Z}_k\}$ converges to the set of equilibria of algorithm \eqref{algo_full}.
Since every equilibrium point of algorithm \eqref{algo_full} is a Lyapunov stable equilibrium followed by Lemma \ref{bound_lemma} (i). It follows from  Lemma \ref{lmm-conv} that $\{\hat{Z}_k\}$ converges to one equilibrium of algorithm \eqref{algo_full}. Equivalently, the sequence $\{\hat{Z}_{i,k}\}$ converges to one same optimal solution for all agent $i$. From Lemma \ref{bound_lemma} (ii), $\sum_{i=1}^m (1-\alpha_i \xi_i) f_i(\bar{Z}_{i,k})=O(\frac{1}{k})$, where $\bar{Z}_{i,k}$ is time average state variable of agent $i$.
$\hfill\square$

\subsection{Proof of Theorem \ref{conv_exp}}
To begin with, the function $f_i(X_i,Y_i)$ can be rewritten as a quadratic form by matrix vectorization. In particular, $vec(ABC)=(A \otimes C') vec(B)$. Then,\\
\begin{align*}
f_{i}=&\frac{1}{2}\left\|Y_i^{ri}\!-\!A_{ri}X_i \right\|_F ^2 \!+\!\frac{1}{2} \left\|Y_{i}A_{ri}'\!-\!X_iE_{li}\!+\!Q_{li} \right\|_F ^2\\
=& \frac{1}{2}\left\| \left[\begin{matrix}E_{ri} & -A_{ri} \end{matrix}\right] \left[ \begin{matrix} Y_i-Y_i^* \\X_i-X_i^*\end{matrix} \right] \right\|_F ^2 \\
 &+ \frac{1}{2} \left\| \left[ \begin{matrix} Y_i-Y_i^* & X_i-X_i^*\end{matrix} \right]   \left[\begin{matrix} A_{ri}' \\ -E_{li}\end{matrix}\right]  \right\|_F ^2\\
=& \frac{1}{2}\left\| vec(\left[\begin{matrix}E_{ri} & -A_{ri} \end{matrix}\right] \left[ \begin{matrix} Y_i-Y_i^* \\X_i-X_i^*\end{matrix} \right]) \right\| ^2 \\
 &+ \frac{1}{2} \left\|vec( \left[ \begin{matrix} Y_i-Y_i^* & X_i-X_i^*\end{matrix} \right]   \left[\begin{matrix} A_{ri}' \\ -E_{li}\end{matrix}\right])  \right\| ^2\\
=& \frac{1}{2}\left\| \left[\begin{matrix}E_{ri} &  -A_{ri} \end{matrix}\right] \otimes E_n vec(\left[ \begin{matrix} Y_i-Y_i^* \\X_i-X_i^*\end{matrix} \right]) \right\| ^2 \\
 &+ \!\frac{1}{2} \! \left\| E_n \!\otimes \!  \left[\begin{matrix}\! A_{ri} & -E_{li}' \!\end{matrix}\right]  \! vec(\left[ \begin{matrix} Y_i\!-\!Y_i^* & X_i\!-\!X_i^*\end{matrix} \right] ) \right\| ^2\\
=& \frac{1}{2}\|C_{i1} vec(\left[ \begin{matrix} Y_i-Y_i^* \\X_i-X_i^*\end{matrix} \right]) \| ^2 \\
 &+ \frac{1}{2} \|C_{i0} vec(\left[ \begin{matrix} Y_i-Y_i^* & X_i-X_i^*\end{matrix} \right] ) \| ^2,
\end{align*}
 where $C_{i1}=\left[\begin{matrix}E_{ri} & -A_{ri} \end{matrix}\right]\otimes E_n  \in \mathbb{R}^{nn_{i}\times 2n^2}$ and $C_{i0}=E_n\otimes \left[\begin{matrix} A_{ri} & -E_{li}'\end{matrix}\right]  \in \mathbb{R}^{nn_{i}\times 2n^2}$. Through finite elementary transformations of $C_{i0}$, we get $f_{i}(X_i,Y_i)= \frac{1}{2}\left\|C_{i1} vec(\left[ \begin{matrix} Y_i-Y_i^* \\X_i-X_i^*\end{matrix} \right]) \right\|_F ^2 + \frac{1}{2} \left\|C_{i2} vec(\left[ \begin{matrix} Y_i-Y_i^* \\ X_i-X_i^*\end{matrix} \right] ) \right\|_F ^2=\frac{1}{2} \left< \mathcal{Z}_{vi}, P_i \mathcal{Z}_{vi}\right>_F$, where $\mathcal{Z}_{vi}=vec\left[ \begin{matrix} Y_i-Y_i^* \\X_i-X_i^*\end{matrix} \right]\in \mathbb{R}^{2n^2}$ and
\begin{align}\label{Pi_def}
P_i=(C_{i1}'C_{i1}+C_{i2}'C_{i2}).
 \end{align}
 In the equation (\ref{Pi_def}), $C_{i2}$ is obtained from $C_{i0}$ through finite elementary transformations.
 \par Define $\mathcal{Z}=\left[\begin{matrix} \mathcal{Z}_{v1}\\ \vdots \\\mathcal{Z}_{vm}\end{matrix}\right]=vec(\hat{Z}-\hat{Z}^*)\in \mathbb{R}^{2mn^2}$. With the formula $f_i(X_i,Y_i)=\frac{1}{2} \left< \mathcal{Z}_{vi}, P_i \mathcal{Z}_{vi}\right>_F$, we rewrite $F(\hat{X},\hat{Y})$ in problem (\ref{opti_problem_origin}) as
 \begin{align}\label{mathbf_P}
  F(\hat{X},\hat{Y})= \frac{1}{2}\left< \mathcal{Z}, \mathbf{P} \mathcal{Z}\right>_F,
  \end{align}
 where $\mathbf{P} = diag\{P_1,\cdots,P_m\} \in \mathbb{R}^{2mn^2\times 2mn^2}$.
 \par Since $\hat{Z}_{k+1}- \hat{Z}^* \!= \!(\hat{Z}_k-\hat{Z}^*)\!-\! \frac{1}{2}\Lambda L_k(\hat{Z}_k-\hat{Z}^*)-\Lambda \nabla_{\hat{Z}_k} F(\hat{Z}_k)$, using the matrix vectorization idea, we obtain
 \begin{align} \label{vec_algo}
 \mathcal{Z}_{k+1}=  \mathcal{Z}_{k}-(\frac{\Lambda  L_k}{2} \otimes E_n )  \mathcal{Z}_{k}-(\Lambda \otimes E_n )  \mathbf{P} \mathcal{Z}_k.
 \end{align}
 where $L_k=\mathcal{L}_k \otimes E_{2n}\in \mathbb{R}^{2mn\times 2mn}$ and $\Lambda =diag\{ \alpha_1\otimes E_{2n}, \cdots, \alpha_m \otimes E_{2n}\} \in \mathbb{R}^{2mn\times 2mn}$. Let $\mathbf{L}_k$ denote $(L_k\otimes E_n)$ and let $\mathbf{\Lambda}$ denote $\Lambda \otimes E_n $. The equation (\ref{vec_algo}) can be further rewritten as
 \begin{align}\label{vec_algorithm}
  \mathcal{Z}_{k+1} = \mathcal{Z}_k- (\frac{\mathbf{\Lambda}\mathbf{L}_k}{2}+\mathbf{\Lambda}\mathbf{P}) \mathcal{Z}_k.
 \end{align}

To proceed, we recall some results on \emph{semistability} of dynamic systems. The notion of \emph{semistability} pertains to a continuum of initial state dependent equilibria, implying that any system trajectory converges to a equilibrium dependent on its initial state \cite{semi_stable}. \emph{Semistability} is the property where every trajectory that starts in a neighborhood of a Lyapunov stable equilibrium converges to a (possibly different) Lyapunov stable equilibrium.
\par A discrete-time switched linear system on $\mathbb{R}^n$ is
\begin{align}\label{sw_lin}
x(k+1)=A(k)x(k),
 \end{align}
 where variable state $x(k)$ evolves by switching over a finite dynamic coefficient matrix sequence $\{A(k)\}$.
 \par By equation (\ref{vec_algorithm}), the proposed algorithm can be considered as a discrete-time switched linear system whose variable is $\hat{Z}_k-\hat{Z}^*$. For the problem (\ref{opti_problem_origin}), if Assumption \ref{assumption1} holds, then the algorithm's corresponding switched linear system has a continuum of equilibria.
\par Based on \emph{semistability} theories of switched linear systems, we get the desired result. Before giving the proof of convergence rate, we present a related lemma, which has been well studied in \cite{semi_switch}.
\begin{lemma} \label{semistability}
	A switched linear system (\ref{sw_lin}) is linearly semistable under arbitrary switching if and only if it is semistable under arbitrary switching.
\end{lemma}
\par In Lemma \ref{semistability}, the linear semistability of a system means that the system's states converge to one stable point at a linear rate.
\par {\bf Proof of Theorem \ref{conv_exp}: }
	Suppose Assumption \ref{connected_ass} holds, $\mathcal{G}_k$ is connected for any iteration time $k$. Then, the variable states converge to one optimal solution under arbitrary switching topologies by Theorem \ref{converge_theo}.
	\par Note that the proposed algorithm can be considered as a discrete-time switched linear system. Because variable states converge to one Lyapunov stable equilibrium of the solution set, the switched linear system is semistable under arbitrary switching topologies. By Lemma \ref{semistability}, the switched linear system corresponding to the proposed algorithm is linearly semistable under arbitrary switching topologies. If Assumption \ref{assumption1} holds, the stable point set of the switched linear system is the set where $\hat{Z}_k=\hat{Z}^*$ and $Z_i^*=Z_j^*$. As a result, we conclude that $\{\hat{Z}_k\}$ converges to one optimal solution $\hat{Z}^*$ at a linear rate.
$\hfill\square$

\section{Simulation}\label{simulation}
In this section, some numerical simulations are presented to show the efficacy of the proposed algorithm.
\par \textit{Example:} For linear time-invariant systems, the Lyapunov test for controllability is $AXA'-X=-BB'$ \cite{linear_theory}. Compared to problem (\ref{target1}), $Q=BB'$, which is symmetric. Suppose that the system information is distributed among undirected graphs containing $5$ agents. The matrices are decomposed as
$A=\left[\begin{matrix} A_{r1} \\ \vdots \\A_{r5}\end{matrix}\right]$ and $Q=\left[\begin{matrix} Q_{l1} \\ \vdots \\Q_{l5}\end{matrix}\right]$.
\par
\begin{figure}
  \centering
  \includegraphics[width=8cm]{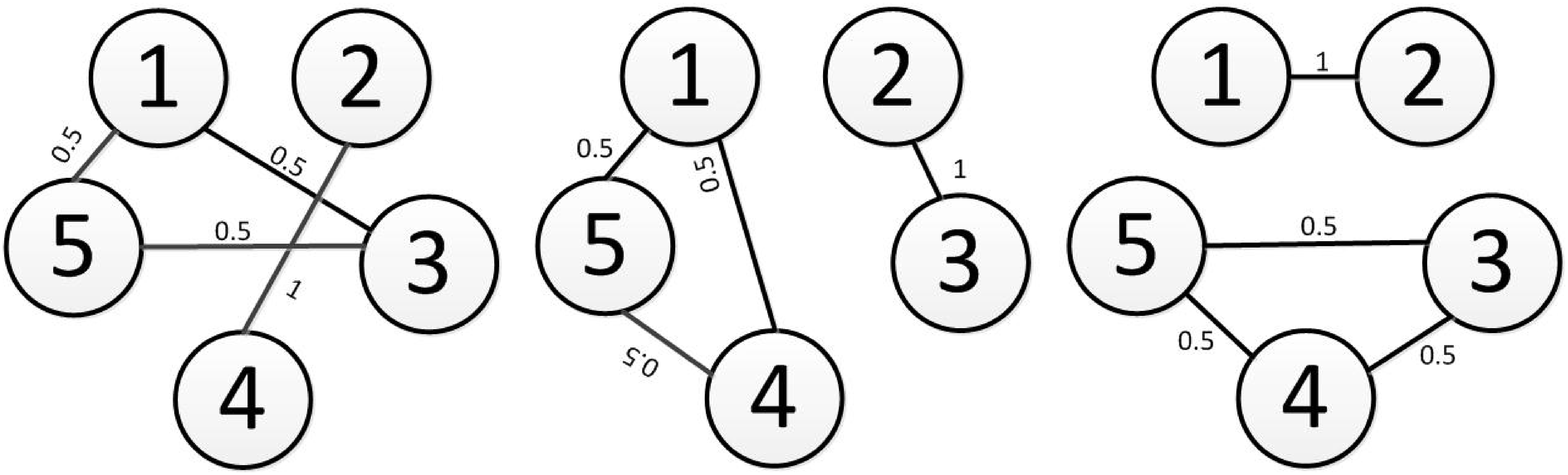}\\
  \caption{uniformly connected graphs}\label{Fig. a}
  \includegraphics[width=8cm]{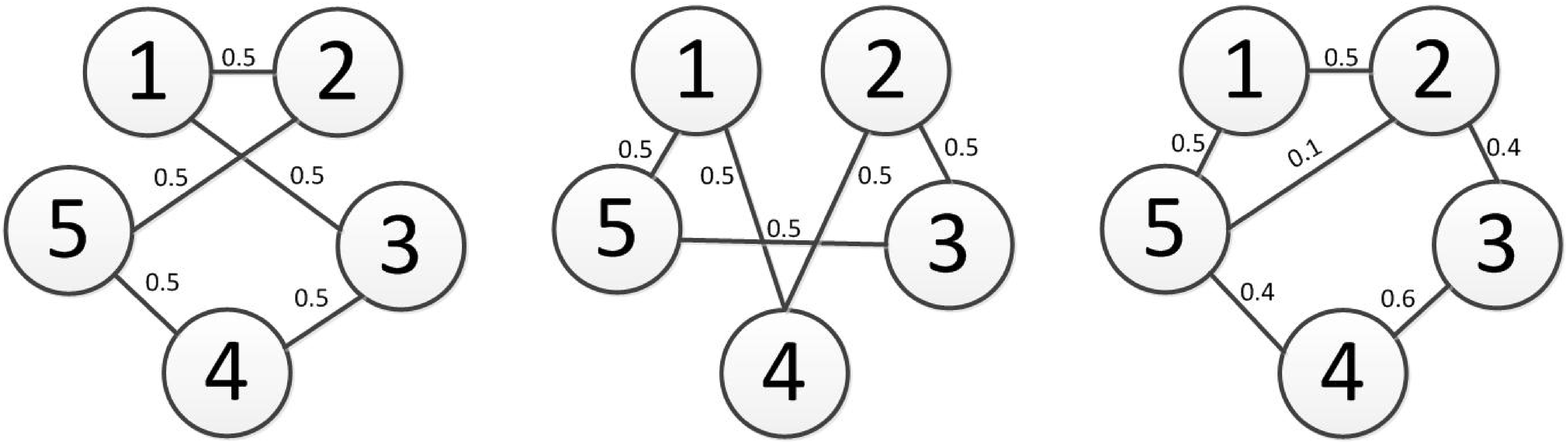}\\
  \caption{connected graphs}\label{Fig. w}
  \centering
\end{figure}

 Every agent $i$ only knows local information of system $A_{ri}$, $Q_{li}$ and communicates with neighbors over undirected graphs. We consider two types of time-varying undirected topologies, which are illustrated in Figs. \ref{Fig. a},\ref{Fig. w}.
For each type of topological graphs, we select one topology randomly from them at every iteration time $k$. Note that there can be many time-varying topologies. Here, we only give three of them for simplicity. Every agent updates its states according to the proposed algorithm (\ref{algo_full}). The constant step sizes for different agents are listed as a vector $sp$.
\par
\begin{figure}
  \centering
  \includegraphics[width=8 cm, height =  6 cm]{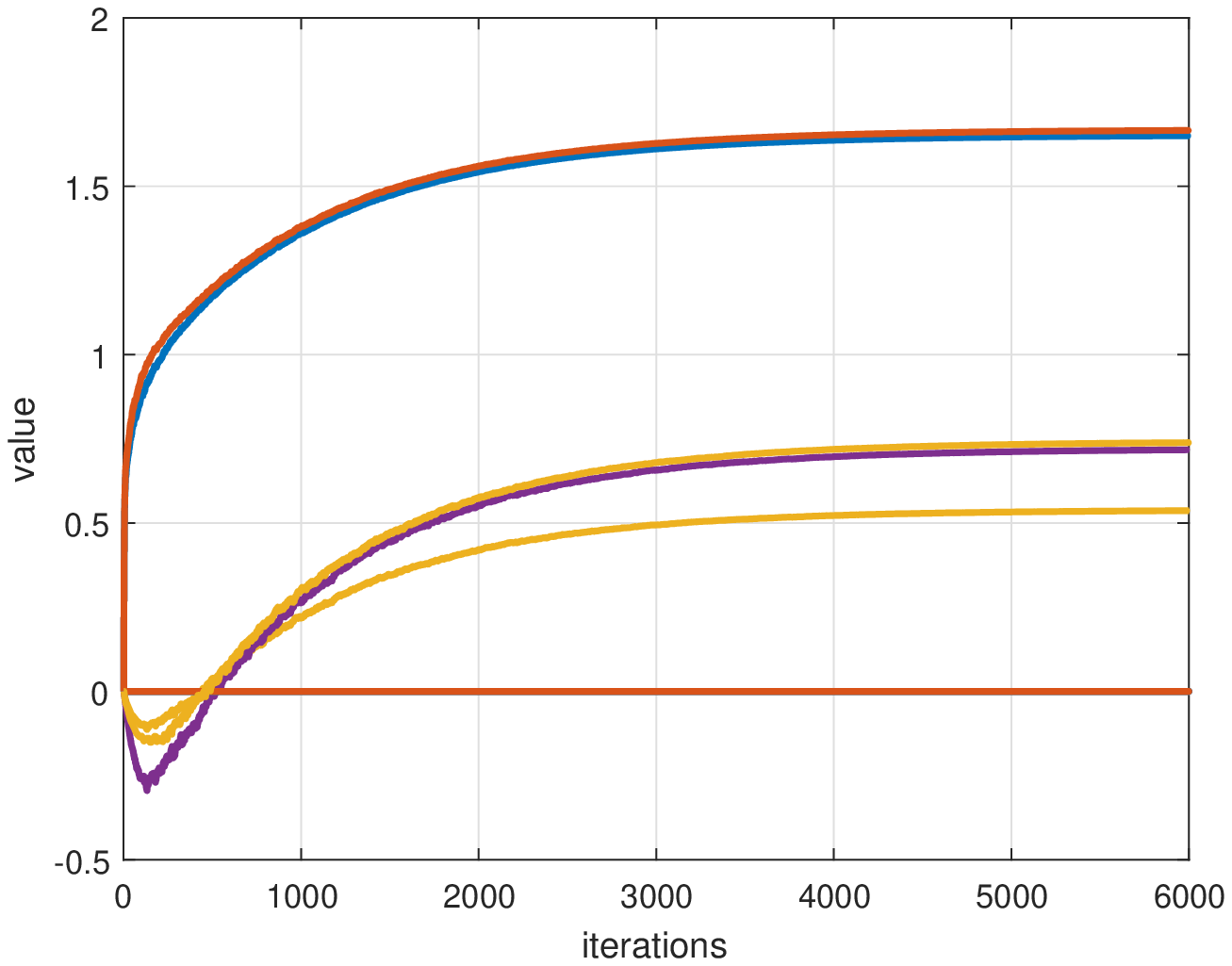}\\
  \caption{The trajectories of row elements of $X$}\label{Fig. 1}
  \includegraphics[width=8 cm, height =  6 cm]{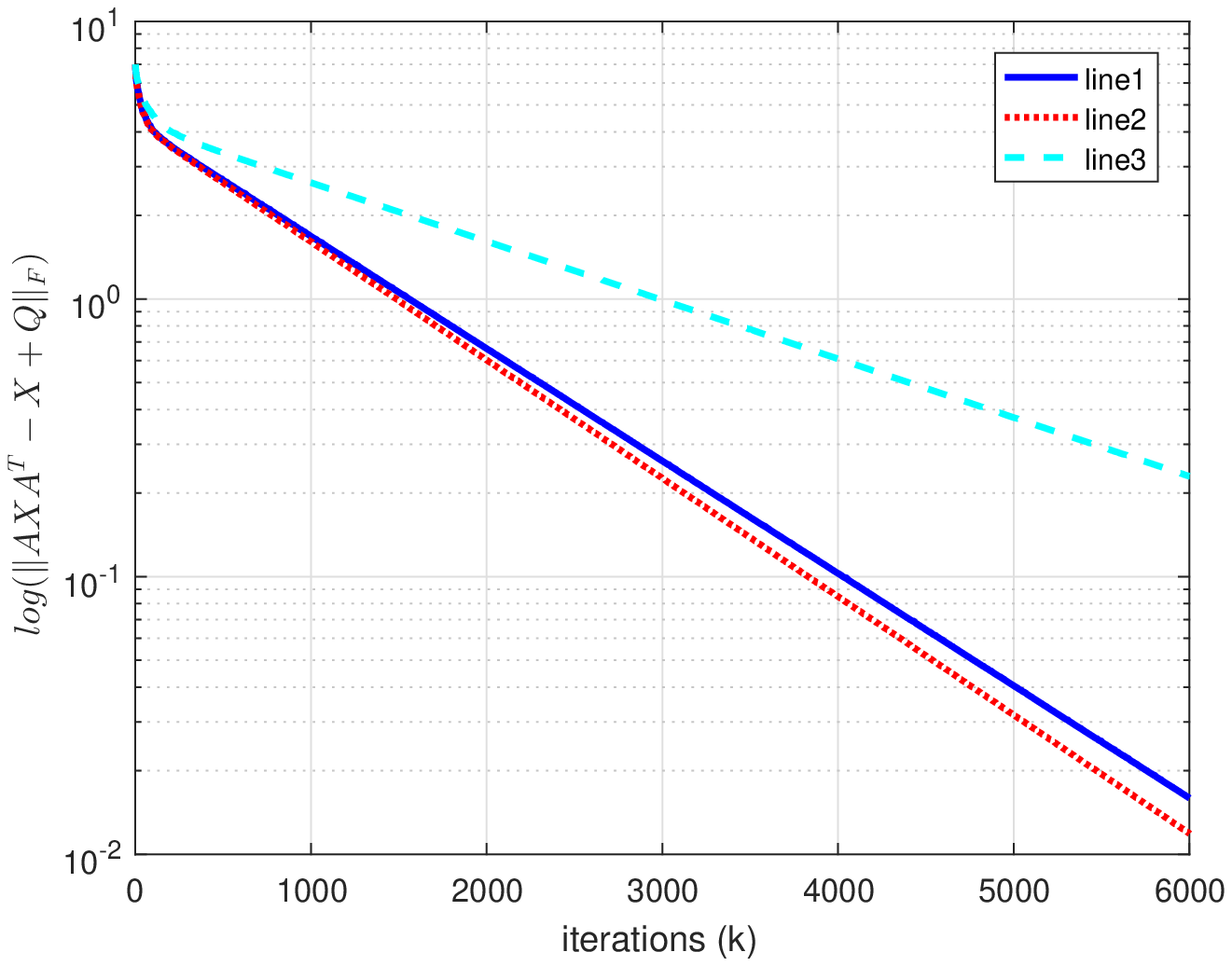}\\
  \caption{The trajectories of $\left\|AXA'-X+Q\right\|_F$ versus iterations}\label{Fig. 5}
  \includegraphics[width=8 cm, height =  6 cm]{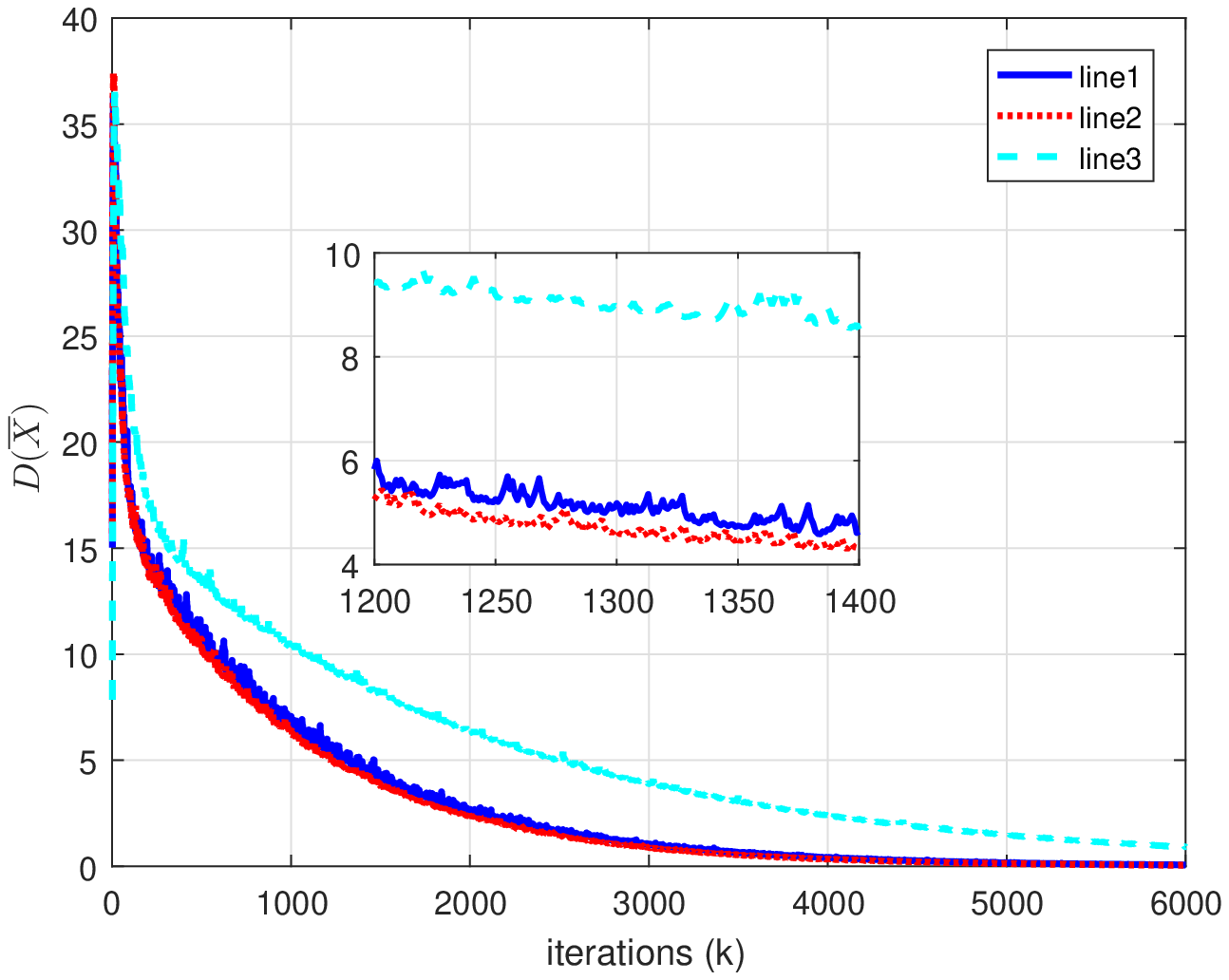}\\
  \caption{The trajectories of $D(\bar{X})$ versus iterations}\label{Fig. 6}
\centering
\end{figure}

 We provide a linear time-invariant system for example, which is controllable by the rank test of the controllable matrix. The coefficient matrices of the selected system $A$ and $B$ are listed in table \ref{matrix_table} and the matrix dimension $n=10$. Using algorithm (\ref{algo_full}) with iteration $k=6000$, we solve the Lyapunov matrix equation to verify the system's controllability. The simulation results are shown in Figs. \ref{Fig. 1}-\ref{Fig. 6}. Furthermore, all eigenvalues of the convergence matrix are positive, which implies that the matrix is positive definite and the linear time-invariant system is controllable. The result is the same as the rank test of controllable matrix, verifying the correctness of the proposed algorithm.
\begin{table*}
        \centering
       \begin{tabular}{|c|c|}
         \hline


         $A=\left[\begin{matrix}0.0061 & 0.1355 & 0.0998 &0.1051 &  0.1085  &  0.0007  &  0.1095  &  0.0817   & 0.0036 &   0.0625 \\0.1492  & 0.1171 & 0.0385 & 0.0708 & 0.0791 & 0.0387 &   0.0523   &0.0311 &  0.0163  &  0.1427 \\ 0.1014  & 0.1360  &  0.0263   & 0.0603 &  0.0087  &  0.0481 &  0.0586  &  0.0165  &  0.1036  &  0.0670\\0.1255 &  0.0109  &  0.0290 & 0.1362&    0.1128 & 0.0086  & 0.0244 & 0.1072  &  0.1397 & 0.0689 \\ 0.1439  & 0.1434  & 0.0378   & 0.0980   & 0.0884  &  0.0410 &  0.1360&  0.1231  &  0.1216  &  0.0058\\ 0.0268  &  0.1180 &   0.1009  &  0.0915  &  0.0303 & 0.1404  &  0.0200 & 0.1167  & 0.0305  &  0.1475\\ 0.0944  &0.1159  &  0.0598  &  0.1209  &  0.1100  & 0.0749 & 0.0108   & 0.1298  &  0.0572  &  0.0684\\    0.1189 & 0.0558  &  0.0391  &  0.1120  &  0.0946 & 0.0952 &  0.1077 &  0.0560  &  0.0840  &0.0900\\ 0.0958  & 0.0067  & 0.0693  &  0.0899 & 0.1106 & 0.0007  &0.0687 & 0.0153 & 0.0499&  0.0388\\0.0828  &  0.0712 &   0.0463 &  0.1241&   0.0535  &  0.1397  &  0.0728  &  0.1038 &   0.1492  &  0.0129\end{matrix}\right]$ &$B=\left[\begin{matrix}1 &1\\1 &0\\0&0\\0 &1\\1 &0\\0& 0\\1 &0\\ 1& 1\\ 0& 1\\ 0& 0 \end{matrix}\right] $\\
         \hline

       \end{tabular}
       \caption{Coefficient matrices A and B}\label{matrix_table}
\end{table*}

\par In Fig. \ref{Fig. 1}, the trajectories of row elements of one agent's estimate $X$ are displayed, which shows that the element values converge to one solution. Fig. \ref{Fig. 5} shows logarithmic trajectories of $\left\|AXA'-X+Q\right\|_F$ over time-varying topologies and Fig. \ref{Fig. 6} shows trajectories of $D(\bar{X})=\sum_{i=1}^m\left\|\sum_{j=1}^m(X_i-X_j)\right\|_F$, $\forall i,j \in \{1,2,3,4,5\}$. The lines in Fig. \ref{Fig. 5} and Fig. \ref{Fig. 6} are demonstrated in Table \ref{table1}. In Fig. \ref{Fig. 5}, all curves over time-varying topologies converge to $0$, implying that the convergent matrix is a solution to the discrete-time Lyapunov matrix equation (\ref{target1}). In addition, the curve with half step sizes $sp/2$ converges slower than the curve over uniformly connected graphs, demonstrating the influence of step sizes on the convergence rate of proposed algorithm. All curves in Fig. \ref{Fig. 6} converge to $0$, which implies that the solution states of all agents converge to the same matrix. Since every graph in connected topologies is connected, the curve converges more smoothly by the detailed diagram in Fig. \ref{Fig. 6}.

\begin{table}
    \centering
    \begin{tabular}{|c|c|c|} 
    \hline
    Lines & Step sizes &Topological graphs\\
    \hline
    line1& $sp$ & uniformly connected\\
    \hline
    line2& $sp$& connected\\
    \hline
    line3& $sp/2$ & uniformly connected\\
    \hline
    \end{tabular}
    \caption{parameter table}\label{table1}
    \centering
\end{table}
\section{Conclusion}\label{conclusion}
The paper has studied the problem of solving large scale discrete-time Lyapunov matrix equations in a distributed way over dynamic multi-agent networks. A distributed algorithm with constant step sizes has been proposed and proofs for convergence properties over time-varying topologies have been presented. Constant step sizes of the proposed algorithm are uncoordinated, depending on every agent's local information. Variable estimates generated by the proposed algorithm converge to an optimal solution of the solution set at a linear convergence rate.

\ifCLASSOPTIONcaptionsoff
  \newpage
\fi



%

\end{document}